\documentclass[11pt]{amsart}
\usepackage{babel}
\usepackage{cjhebrew}

\usepackage[dvipsnames]{xcolor}
\usepackage{amssymb,verbatim}
\usepackage{amsmath,amsfonts,enumitem,bm}
\usepackage[mathscr]{euscript} 
\usepackage{amsthm}
\usepackage{url}
\usepackage{graphicx} 
\usepackage{enumitem} 
\usepackage{hyperref} 
\hypersetup{colorlinks} 
\usepackage{bbm} 
\usepackage{bm} 
\usepackage{mathrsfs} 
\usepackage{cancel} 
\usepackage{ytableau} 
\usepackage{mathtools}

\usepackage[colorinlistoftodos]{todonotes}

\RequirePackage{cleveref}
\usepackage{hypcap}
\hypersetup{colorlinks=true, citecolor=darkblue, linkcolor=darkblue}
\definecolor{darkblue}{rgb}{0.0,0,0.7}
\newcommand{\darkblue}{\color{darkblue}}

\definecolor{darkred}{rgb}{0.68,0,0}
\newcommand{\darkred}{\color{darkred}}

\definecolor{darkgreen}{rgb}{0,.38,0}
\newcommand{\darkgreen}{\color{darkgreen}}

\definecolor{magenta}{rgb}{.51, 0, .51}
\newcommand{\magenta}{\color{magenta}}

\newcommand{\defn}[1]{\emph{\darkblue #1}}
\newcommand{\defna}[1]{\emph{\darkred #1}}

\newcommand{\defng}[1]{\emph{\darkgreen #1}}
\newcommand{\defnm}[1]{\emph{\magenta #1}}

\setlist[enumerate]{
	label=\textnormal{({\roman*})},
	ref={\roman*}}

\makeatletter
\def\th@plain{%
	\thm@notefont{}
	\itshape 
}
\def\th@definition{%
	\thm@notefont{}
	\normalfont 
}
\makeatother

\newtheorem{thm}{Theorem}[section]

\newtheorem*{claim*}{Claim}

\newtheorem{conv}[thm]{Convention}

\theoremstyle{definition}

\numberwithin{figure}{section}
\numberwithin{equation}{section}

\def\nn{\mathbb N}
\def\cc{\mathbb C}

\def\la{\lambda}

\def\al{\alpha}

\def\<{\langle}
\def\>{\rangle}

\def\0{{\mathbf 0}}

\def\.{\hskip.06cm}
\def\ts{\hskip.03cm}

\def\pt{\partial}

\def\bx{{\textbf{\textit{x}}}}

\def\.{\hskip.06cm}
\def\ts{\hskip.03cm}

\newcommand{\textsu}[1]{\textup{\textsf{#1}}}

\newcommand{\ComCla}[1]{\textup{\textsu{#1}}}

\newcommand{\sharpP}{\ComCla{\#P}}
\newcommand{\SP}{\ComCla{\#P}}

\newcommand{\NP}{\ComCla{NP}}

\newcommand{\BPP}{\ComCla{BPP}}
\newcommand{\BQP}{\ComCla{BQP}}

\newcommand{\coNP}{\ComCla{coNP}}
\newcommand{\E}{\ComCla{E}}
\renewcommand{\P}{\ComCla{P}}

\newcommand{\PH}{\ComCla{PH}}

\newcommand{\AM}{\ComCla{AM}}

\newcommand{\NE}{\ComCla{NE}}
\newcommand{\coNE}{\ComCla{coNE}}

\def\SP{\sharpP}

\def\GRH{\textup{\sc GRH}}
\def\ERH{\textup{\sc ERH}}

\def\HNP{\textup{\sc HNP}}

\def\MVA{\textup{\sc MVA}}
\def\IWA{\textup{\sc IWA}}
\def\ETH{\textup{\sc ETH}}

\def\poly{{\P}}

\newcommand{\Sch}{\mathfrak{S}} 

\begin{document}

\title[Positivity of Schubert Coefficients]{Positivity of Schubert Coefficients}

\author[Igor Pak \. \and \. Colleen Robichaux]{Igor Pak$^{\ts \bigstar\ts \blacktriangle}$  \. \and \.  Colleen Robichaux$^{\ts\bigstar\ts \blacklozenge}$ }

\thanks{\thinspace ${\hspace{-.45ex}}^\bigstar$ Department of Mathematics,
UCLA, Los Angeles, CA 90095, USA. Email:~\texttt{\{pak,robichaux\}@math.ucla.edu}}


\thanks{$^\blacktriangle$ Supported by the NSF grant CCF-2007891.}
\thanks{$^\blacklozenge$ Supported by the NSF MSPRF grant DMS-2302279.}
\thanks{\today}

\begin{abstract}
\emph{Schubert coefficients} \ts $c^w_{u,v}$ \ts are structure constants
describing multiplication of Schubert polynomials.  Deciding positivity of
Schubert coefficients is a major open problem in Algebraic Combinatorics.
We prove a positive rule for this problem based on two standard assumptions.
\end{abstract}

\maketitle

\begin{cjhebrew}
 drK 'mwnh b.hrty
\end{cjhebrew}


\section{Introduction}\label{s:intro}

In the era of specialization we live in, it is common for
experts in one area of mathematics to be unable to understand
and evaluate results in another area.  Occasionally, things turn
for the worse, when one is unable to understand the \emph{definitions},
or even the \emph{notation}.  The results then become expressions
in a foreign language that take time, effort and dedication to translate.

Most recently, we found ourselves in an unenviable position of having
made progress on a major open problem in Algebraic Combinatorics
\cite{PR-HN}, yet our results are stated in the language
of Computational Complexity that was deep and technical enough
to be inaccessible to most people interested in the problem.
Given a choice of being understood but unmotivated vs.\ being
motivated and misunderstood, we chose the former.

The main goal of our study is to give a {\bf \defnm{combinatorial
interpretation for Schubert coefficients}}, arguably the most important
open problem in Schubert calculus.  In \cite{PR-HN}, we analyze the
{\bf \defnm{Schubert positivity problem}}, traditionally a stepping
stone towards a combinatorial interpretation and a major open problem in
own right.  The main result in \cite{PR-HN} is an upper bound for the
complexity class of the Schubert positivity problem.  We presented the
results in a conventional CS theory style, while the proofs used some
technical algebro-geometric arguments.

In this note, we address the Algebraic Combinatorics audience,
giving implications of our main result in \cite{PR-HN} using
conventional combinatorics language.  More precisely, we show
that the Schubert positivity problem has a {\bf \defnm{positive rule}}
\ts assuming two standard hypotheses.
%
To streamline the presentation, we omit most of the
Schubert calculus background  
that was included in \cite{PR-HN},
as well as the extensive complexity background which
can be found, e.g., in \cite{Gol08,Wig19}.
%

\medskip

\section{The results}\label{s:main}

\subsection{Background} \label{ss:main-def}
Below we give a brief review of the necessary background.
We refer to~\cite{Mac91,Man01} for introductory surveys, to \cite{Knu22}
for an overview of recent results, to \cite{AF} for geometric aspects,
and to \cite[$\S$10]{Pak-OPAC} for computational complexity aspects.

\defng{Schubert polynomials} \ts $\Sch_w\in \nn[x_1,\ldots,x_n]$ \ts indexed by
permutations $w \in S_n\ts$. They are celebrated generalizations of \emph{Schur polynomials}
but are not symmetric in general.  Schubert polynomials were introduced by
Lascoux and Sch\"{u}tzenberger in~1982 to represent
cohomology classes of Schubert varieties in the complete flag variety,
and have been intensely studied from algebraic, combinatorial, geometric,
and (more recently) complexity points of view.

Let \. $\Sch_{w_\circ} : = \. x_1^{n-1} x_2^{n-2} \. \cdots  \. x_{n-1}\ts$,
{where}  \ts $w_\circ := (n,n-1,\ldots,1)\in S_n$ \ts is the \emph{long permutation}.
%
\defn{Schubert polynomials} are defined recursively: \.
$\Sch_{ws_i} := \pt_i \ts \Sch_w$  \. for all \emph{descents} \ts $w(i)> w(i+1)$.
Here~$\ts \pt_i$ \ts is the \emph{divided difference operator} \ts defined as \.
$$\pt_i F \, := \, \frac{F - s_i F }{x_i -x_{i+1}}\,, \quad \text{where} \quad
s_i \. := \. (i,i+1) \ts \in \ts S_n
$$
is the transposition which acts on \ts
$F\in \cc[x_1,\ldots,x_n]$ \ts by transposing variables. Schubert polynomials coincide
with Schur polynomials for 
permutations with at most one descent.

It is known that Schubert polynomials have integral coefficients:
\ts $[\bx^\al] \Sch_w \in \nn$ \ts called \emph{Schubert--Kostka numbers},
which generalize the usual \emph{Kostka numbers}.
Schubert--Kostka numbers have combina-torial interpretations in terms of
\emph{pipe dreams} (\emph{RC-graphs}) and \emph{bumpless pipe dreams}.

Schubert polynomials
form a basis in the ring of polynomials, with integral structure constants:
$$
\Sch_u \cdot \Sch_v \, = \, \sum_{w \ts \in \ts S_\infty} \. c^w_{u,v} \. \Sch_w\ts.
$$
It is known that \ts $c^w_{u,v} \in \nn$ \ts for all \ts $u,v,w\in S_\infty$\ts,
as they have both geometric and algebraic meanings, which generalize
the number of intersection points of lines, see e.g.\ \cite{AF,BSY}.  These integers are
called \defn{Schubert coefficients} and generalize the \emph{Littlewood--Richardson} (LR-)
\emph{coefficients}.

LR-coefficients have over 20 combinatorial interpretations, see \cite[$\S$11.4]{Pak-OPAC}.
By comparison, Schubert coefficients have two different \emph{signed
combinatorial interpretations}, see \cite[Prop.~10.2]{Pak-OPAC} and \cite{PR24}.
Whether Schubert coefficients have a (unsigned) combinatorial interpretation is a major open
problem, see \cite[Problem~11]{Sta00} and \cite[$\S$1.4]{Knu22}.

\defn{Schubert positivity} \ts is a problem whether Schubert coefficients
are positive:
\begin{equation}\label{eq:Schu-pos}
\big\{c^w_{u,v} \. >^? 0 \big\}.
\end{equation}

%
This problem was heavily studied and resolved in special cases.
Notably, Knutson \cite{Knutson01}, Purbhoo \cite{Purbhoo06}, Billey and Vakil \cite{BV08},
and, most recently, St.~Dizier and Yong \cite{StDY22}, gave necessary 
conditions for Schubert positivity.

\smallskip

\subsection{Positive rule} \label{ss:main-comb}
In Algebraic Combinatorics, a {\bf \defna{combinatorial interpretation}} \ts is often informally
defined as
\begin{equation}\label{eq:ast}
\text{\emph{``manifestly nonnegative formula''}, \ts see e.g. \cite[$\S$1.4]{Knu22}. }
\end{equation}
This notion is 
best understood in context, especially
when compared to standard combinatorial interpretations of Kostka numbers and LR-coefficients.

First, recall that a \emph{closed formula} \ts is best defined as a function
computable by an efficient algorithm, see \cite{Wilf82} (cf.~\cite[$\S$1.1]{Sta99}).
Now, the word ``formula''
in \eqref{eq:ast} is meant to be a summation over combinatorial objects of
closed formulas, rather than a stand alone closed formula.
In fact, neither Kostka numbers nor LR-coefficients can be computed efficiently,
see a discussion in \cite[$\S$5.2]{Pan23}.

Second, the word ``manifestly'' in this context is a reference to
a sum of positive terms, since a signed summation is not self-evidently positive.
Taken together, \eqref{eq:ast} can be rephrased as
%
\begin{equation}\label{eq:ast-ast}
\text{\emph{``summation over combinatorial objects of positive closed formulas''}. }
\end{equation}

In \cite[$\S$11.4]{Pak-OPAC}, the first author clarified this definition further
as membership in $\SP$, a standard computational complexity class.    The meaning of
``combinatorial interpretation'' is then:
\begin{equation}\label{eq:circ}
\text{\emph{``the number of combinatorial objects, which are poly-time verifiable''}. }
\end{equation}
Compared to $\eqref{eq:ast-ast}$, this specifies coefficients of the summations to be~$1$,
and restricts the type of combinatorial objects to those whose validity can be
verified in polynomial time.  Still, definition $\eqref{eq:circ}$ is broad enough to
include the numerous combinatorial interpretations of LR-coefficients mentioned above.
For example, one can verify in polynomial time if a given combinatorial object is indeed a
\emph{LR-tableaux}, a \emph{Gelfand--Tsetlin pattern},
a \emph{Berenstein--Zelevinsky triangle}, a \emph{Knutson--Tao puzzle}, etc.

The \emph{positivity} of combinatorial numbers is often
viewed in a similar spirit.  Occasionally, the positivity can be decided in polynomial
time; this is the case of Kostka and LR-coefficients (see e.g.\ \cite[$\S$5.2]{Pan23}).
Other times, this is known to be impossible,  e.g.\ the positivity problem
for Kronecker coefficients is $\NP$-hard for Kronecker coefficients (ibid.)
In many such cases, a {\bf \defna{positive rule}} \ts is used to establish the positivity:
%
\begin{equation}\label{eq:di}
\text{\emph{``positivity is given by a combinatorial object, which is poly-time verifiable''}. }
\end{equation}

In other words, if the desired number is positive, then there is a combinatorial
object as above, and vice versa.
In particular, if the positivity is decidable in polynomial time, this algorithm by
itself gives a positive rule (with input as the combinatorial object).
In the language of Computational Complexity, this says that the problem is in~$\NP$.
Clearly,
\emph{every} combinatorial object counted by a
combinatorial interpretation $\eqref{eq:circ}$ gives a positive rule $\eqref{eq:di}$.
Thus a positive rule can in principle
be easier to establish than a combinatorial interpretation.

For example, a single $3$-coloring of a graph is a positive rule for 
of the number of $3$-colorings.  Similarly, a single LR-tableau is a combinatorial
rule for positivity of a LR-coefficient.  Let us emphasize that the existence of
a such tableau is both a necessary and sufficient condition.  On the other hand,
following \cite[$\S$3.2]{Purbhoo06}, a single \emph{winning root game} suffices to
show that a given Schubert coefficient is positive.  This is only a sufficient
condition for positivity, thus not a positive rule.  Similarly, the existence
of a \emph{permutation array} is a necessary condition for positivity of Schubert
coefficients \cite{BV08}, thus not a positive rule again.

We conclude by noting that there are nonnegative combinatorial functions with
no positive rules (unless polynomial hierarchy \ts $\PH$ \ts collapses),
see \cite{CP-AF,IP22}.  Notably, in \cite{IPP24} we showed that non-vanishing problem
of the $S_n$ character \. $\big\{\ts |\chi^\la(\mu)|>^?0\big\}$ \ts does not have a
positive rule, implying that the squared character does not have a
combinatorial interpretation.

\smallskip

\subsection{Main theorem}\label{ss:main-thm}
Our main result is a positive rule for 
Schubert coefficients under two assumptions.  This is the first general claim
towards the long sought combinatorial interpretation.

\smallskip

\begin{thm}
\label{t:main}
Assuming  \ts $\GRH$ \ts and \ts $\MVA$,  the Schubert positivity problem
 \eqref{eq:Schu-pos}
has a positive rule.
\end{thm}

\smallskip

Here the $\GRH$ is the \defn{Generalized Riemann Hypothesis} \ts which states that
all nontrivial zeros of $L$-functions \ts $L(s,\chi_k)$ \ts have real part~$\frac12$.
In fact, the \emph{Extended Riemann Hypothesis} ($\ERH$, see e.g.\ \cite[$\S$6]{BCRW08}),
or an even weaker assumption in \cite[Thm~2(2)]{Roj03}, also suffice.

The $\MVA$ is the \defn{Miltersen--Vinodchandran Assumption}, a strong {derandomization
assumption} \ts which implies a collapse of two complexity classes: \ts $\NP = \AM$.
Formally, $\MVA$ \ts states that some language in \ts $\NE{}\. \cap \. {}\coNE$ \ts requires
nondeterministic circuits of size \ts $2^{\Omega(n)}$.
This assumption was introduced in \cite[Thm~1.5]{MV05} as an interactive proof
analogue of the \emph{Impagliazzo--Wigderson Assumption} ($\IWA$), that
some problem in \ts $\E$ \ts (say, {\sc SAT}), requires circuits of size \ts $2^{\Omega(n)}$.
The $\IWA$ is a classic derandomization assumption which implies \ts $\P= \BPP$, i.e.,
that all probabilistic polynomial time algorithms can be made deterministic.

The assumption \ts $\MVA$ \ts is best viewed as a substantial
strengthening of the \ts $\P\ne \NP$ \ts conjecture, far beyond the
\emph{Exponential Time Hypothesis} ($\ETH$).  We refer to \cite[$\S$8.3, $\S$9.1]{Gol08}
and \cite[$\S$7.2]{Wig19} for more of these results in the context of computational complexity,
and to \cite[$\S$7]{MV05} for prior work, stronger assumptions, and an overview of followup results.

\smallskip

\begin{proof}
We deduce the result from \cite{PR-HN} and \cite{MV05}.  We showed in \cite[Lemma~1.10]{PR-HN},
that a modification of the \emph{lifted formulation} given in \cite{HS17} can be used to show that the
Schubert positivity problem reduces to \emph{Parametric Hilbert's Nullstellensatz} ($\HNP$).
Here $\HNP$ asks if the polynomial system \. $f_1 = \. \ldots = f_m = 0$ \. has a solution over
\ts $\overline{\mathbb{C}(y_1,\ldots,y_k)}$, where \ts $f_i \in \mathbb{Z}(y_1,\ldots,y_k)[x_1,\dots,x_s]$
\ts for all \ts $1\le i \le m$.
In \cite[Thm~1]{A+24}, the authors extend Koiran's celebrated result \cite{Koiran96}, to show that
$\HNP$ is in $\AM$ assuming $\GRH$.  In \cite[Thm~1.5]{MV05}, the authors 
show that \ts $\AM=\NP$ \ts assuming $\MVA$.  Thus, the Schubert positivity problem \eqref{eq:Schu-pos}
is in \ts $\NP$ \ts given both assumptions, as desired.
\end{proof}

\medskip

\section{Discussion}\label{s:disc}

\subsection{The meaning of Theorem~\ref{t:main}} \label{ss:disc-mean}
This is the first positive rule for Schubert positivity in full generality.  It is also the first general result in favor of
Schubert coefficients having a combinatorial interpretation, which would contradict
Conjecture~10.1 in \cite{Pak-OPAC}.  The theorem by itself gives no indication
as to whether Schubert positivity is $\NP$-hard asked in \cite[Question~4.3]{ARY19},
but if the answer is positive (as we expect), this suggests that \eqref{eq:Schu-pos} is
$\NP$-complete.\footnote{Here and elsewhere we are assuming
that permutations are given in their natural presentation. }

We deduce Theorem~\ref{t:main} from a known complexity theoretic result \cite{MV05},
a recent breakthrough \cite{A+24}, and one of our results in \cite{PR-HN}.  As we
explain in \cite{PR-HN}, our result extends verbatim to root systems $B$ and $C$,
but not to~$D$, since the corresponding lifted formulation in the latter case
does not satisfy $\HNP$'s requirements.

Whether the positive rule we obtain in Theorem~\ref{t:main} is especially
\emph{combinatorial} is in the eye of the beholder.  Roughly, the positive rule one
gets from our proof combined with the algorithm in \cite{A+24,Koiran96},
consists of solutions of the lifted formulation
system over multiple primes whose existence is guaranteed by $\GRH$.  The variables
happen to be matrix entries of matrices which can be interpreted as flags over the
finite fields.  Now, since one needs random bits to test the resulting polynomial
identities, these are provided by a \emph{pseudorandom number generator}
constructed in \cite{MV05} using combinatorial tools from hard problems such
as {\sc SAT} or {\sc HamiltonCycle}, which are believed to have exponential circuit
complexity.  All in all, it's an involved rule that is deeply combinatorial in its
construction.

The use of the $\GRH$ assumption is somewhat unfortunate, and possibly
avoidable by a more involved  number theoretic argument.
The use of the $\MVA$ assumption may seem surprising but should not,
since these complexity assumptions are
often left unstated.  Indeed, given a \emph{signed} combinatorial interpretation for
Schubert coefficients, \emph{every} combinatorial interpretation as in \eqref{eq:ast}
becomes trivial from the complexity point of view, if \ts $\poly=\NP$, for example.
For the same reason one can never prove a ``nonexistence of a combinatorial
interpretation'' result without some kind of complexity assumptions.

Finally, one can ask if our rule for Schubert positivity \eqref{eq:di}
should be viewed as an indication towards a combinatorial interpretation
\eqref{eq:circ}.  While we favor a positive answer, there are two negative
arguments to keep in mind.  First, note that Koiran's algorithm
\cite{Koiran96} at the heart of our construction is based on Hilbert's
Nullstellensatz and does not extend to counting the number of solutions,
so a completely different approach is needed.

Second, we now know of a
natural problem, the \emph{defect of the strong Mason's inequality} \ts
for binary matroids, see \cite[$\S$14]{CP-SY}, whose the positivity
is $\NP$-complete (see \cite[Cor.~15.3]{CP-SY}), while the counting
is (conjecturally) not in~$\SP$.  Since Schubert positivity is
(conjecturally) $\NP$-hard (see \cite[Conj.~1.5]{PR-HN}), this
argument suggests that the problem of combinatorial interpretations
is more intricate than we initially thought.

\smallskip

\subsection{The meaning of derandomization} \label{ss:disc-derandom}
Note that $\GRH$ is probably the most famous conjecture in all of mathematics,
with countless references extolling its powers (see e.g.\ \cite{BCRW08}),
and there is a universal belief that it holds (cf.\ \cite{Far22}).
By contrast, $\MVA$ and the whole area of \defng{Derandomization} \ts
may seem unfamiliar and even counterintuitive.  While we cannot give
it justice in a few paragraphs, let us try nonetheless to give some motivation
for this direction.  For a proper introduction we refer to lecture notes
\cite{Vad11}, textbook \cite{Gol08}, and a book-length survey~\cite{Wig19}.

On a basic level, derandomization aims to ``simulate random bits'' by an
algorithm which produces a binary sequence that can be used by a given
probabilistic algorithm as if it was truly random.
While the idea is rather old and practical, theoretical results
are both difficult and technical.
In a major breakthrough, Nisan--Wigderson (1994) and later
Impagliazzo--Wigderson (1997), gave reasonable hardness assumptions
to imply that \ts $\BPP=\poly$.
Later work extended this \defng{hardness-randomness tradeoff} \ts to
other points on the hardness spectrum, and to other complexity classes.

Most relevant to this work, Miltersen--Vinodchandran \cite{MV05}
generalized this phenomenon to interactive proof systems (of which $\AM$
is the first nontrivial class).   They extended prior
work by Arvind--K\"obler (2001) and Klivans--van~Melkebeek (2002),
which proved the \ts $\AM=\NP$ \ts collapse under stronger assumptions.
In the opposite direction, the authors weakened \ts $\MVA$ \ts at
the expense of making a smaller collapse from $\AM$ to a
quasipolynomial version of~$\NP$, see \cite[Thm~1.9]{MV05}.

Part of the effort to derandomize $\AM$ comes from applications
to {\sc GraphIsomorphism}, after it was shown by
Goldreich--Micali--Wigderson (1986), that {\sc GraphNonIsomorphism}
is in~$\AM$.  Since {\sc GraphIsomorphism} is trivially in~$\NP$,
the $\MVA$ implies that this problem is in \ts $\NP \ts \cap \ts \coNP$,
the class which contains {\sc IntegerFactoring}.  Babai's unconditional
quasipolynomial time algorithm (2016) was a major development
that opened a possibility of a poly-time algorithm.

A long chain of successes in the area of derandomization led to a large
$98\%$ majority belief in the $\BPP=\P$ conjecture,
according to the latest poll of experts \cite[p.~21]{Gas19}.
This is shy of $99\%$ majority that $\P\ne \NP$ but very close to
a universal belief.  By contrast, only $70\%$ of experts believe that
{\sc GraphIsomorphism} is in $\P$, and only $31\%$ of experts
believe that {\sc IntegerFactoring} is in~$\P$.  This suggests serious
doubts in the \ts $\NP \cap \coNP =\P$ \ts conjecture, and a substantial
belief that quantum algorithms are more powerful that deterministic
algorithms: \ts $\BQP\ne \P$ \ts (ibid.)

If there is one conclusion to be made from these unscientific polls,
it's that there is a widespread belief among complexity theorists
that probabilistic algorithms are in fact no more powerful than
deterministic algorithms, even if this might take centuries to prove
(since \ts $\MVA$, \ts $\IWA$ \ts and other derandomization assumptions
are stronger than \ts $\P\ne \NP$).  This may feel contrary to the obvious
and numerous successes of Monte Carlo algorithms in many applications,
but suggests we should take these derandomization assumptions
very seriously. 



\medskip

\section{Two critiques} \label{s:critiques}

\subsection{Not interesting}   \label{s:critiques-int}
First, one can argue that the positive rule for Schubert positivity
given by the proof above is substantially different from any combinatorial
objects that had been studied before, that such notion is not useful
in applications and thus not interesting.  While this argument can
be neither formalized nor refuted, it is still worth addressing for
the argument's sake.

The first author argued in \cite{Pak-OPAC} that the complexity class $\SP$
is the right notion to define a ``combinatorial interpretation'', as
encompassing and generalizing all standard examples in the area. We emphasized
that without a proper definition one can never argue \emph{against}
existence of combinatorial interpretation.  We also noted that this
``not in $\SP$'' approach had already been fruitful in several
interesting cases \cite{CP-AF,CP-coinc,IP22,IPP24}.  As we then explain
in \cite[$\S$2.2{\small (2)}]{PR24}, our main result can also be viewed
as the opposite direction, as an obstacle to having a proof of this kind.

Curiously, the literature rarely addresses the possibility of a
positive but extremely involved solution.  One exception is the problem
whether a given knot is an \emph{unknot}, see e.g.\ \cite[$\S$4.6]{Pak-OPAC}.
In this case, a positive rule ($\NP$ witness)
is a sequence of Reidemeister moves which can be easily verified,
but the crux of the argument is that such sequence always exists of
polynomial length.  Even more impressively, testing whether a knot is
{\bf \defng{not}} \ts an \ts \emph{unknot}, is also in~$\NP$.  This is proved by
a highly sophisticated argument using Gabai's construction of taut foliations
as a way of certifying the genus of a knot, see a discussion in \cite[$\S$3]{Lac17}.
Viewed through the prism of combinatorial interpretations, this construction
is highly unintuitive compared to Reidemeister sequences, and at least as technical
as our construction.

As we mentioned above, in the area of Enumerative Combinatorics,
there is an old tradition of using involved efficient algorithms for a
 ``closed formula'' \cite{Wilf82}.  See also \cite{Pak18}
for an extensive overview of this approach.  Once you cross
the bridge from ``nice closed formulas'' to poly-time algorithms,
it is not a big leap to go from ``nice combinatorial interpretations''
\eqref{eq:ast}  to poly-time certificates~$\eqref{eq:circ}$.
Arguably, this is how one turns art into science.

Finally, the relative lack of progress in the positive direction and the
exceedingly cumbersome nature of recent combinatorial interpretations
for certain families of Schubert coefficients \cite{KZ17,KZ23}, suggest
that it is unlikely there is a nice and fully satisfactory combinatorial
interpretation, or even a positive rule.  If history is the guide,
one might want to let go of this dream and consider our more general notions:

\medskip

\begin{center}\begin{minipage}{13.1cm}%
\emph{``Bien des fois déjà on a cru avoir résolu tous les
problèmes, ou, tout au moins, avoir fait l'inventaire de
ceux qui comportent une solution. Et puis le sens du mot
solution s'est élargi, les problèmes insolubles sont
devenus les plus intéressants de tous et d'autres
problèmes se sont posés auxquels on n'avait pas songé.''
\ts {\rm (Henri~Poincar\'{e}, 1913)}}
\end{minipage}\end{center}

\smallskip

\begin{center}\begin{minipage}{13.1cm}%
``Many times already men have thought that
they had solved all the problems, or at least that
they had made an inventory of all that admit of
solution. And then the meaning of the word solution
has been extended; the insoluble problems have
become the most interesting of all, and other problems
hitherto undreamed of have presented themselves.'' \ts {\rm \cite[p.~22]{Poi14}}
\end{minipage}\end{center}

\smallskip

\subsection{Lack of certainty} \label{s:critiques-cert}
One can also argue that our results are quite weak since they rely on two
major unproven hypotheses, each of them stronger than a \ts
{\tt Millennium Problem}.  This is true, of course, but somewhat misleading
since in other areas assuming conjectures is part of the culture.
Notably, there are hundreds of papers in Number Theory and
its applications which assume variants of the GRH as well as other
standard conjectures:  the \emph{ABC conjecture}, the
\emph{Bateman--Horn conjecture}, the
\emph{Birch and Swinnerton-Dyer conjecture},
the \emph{Cohen--Lenstra heuristics},
\emph{Schanuel's conjecture}, the
\emph{Shafarevich--Tate conjecture},
etc.  Sometimes, later arguments
manage to remove or weaken these conjectures as an assumption,\footnote{For example,
this is what happened to the ``not an unknot'' result discussed above,
as the original Kuperberg's argument used the GRH in a manner similar
to \cite{Koiran96}.  Later, Agol removed this assumption, see \cite[$\S$3]{Lac17}.}
but all such results are still considered very valuable.

Similarly, in Computation Complexity, Cryptography and related areas of
Theoretical Computer Science, a large majority of results use \emph{some}
complexity assumptions.  In addition to \ts $\poly\ne\NP$, \ts $\BPP = \poly$ \ts and
non-collapse of the polynomial hierarchy~$\PH$,  
standard assumptions include the \emph{decisional Diffie--Hellman assumption},
the \emph{existence of one-way functions assumption},
the \emph{exponential time hypothesis},
the \emph{learning with errors assumption},
the \emph{small set expansion hypothesis},
the \emph{unique games conjecture}, etc.   Despite this apparent lack of certainty,
the area learned to survive and prosper decades ago.

Returning to Schubert positivity, it is certainly possible that there is an
unconditional positive rule of explicit combinatorial nature.  Unfortunately,
the known special cases are most definitely not general enough to remain
hopeful.  Still, we are completely persuaded by Theorem~\ref{t:main} that there
is \emph{some} \ts (unconditional) positive rule for this problem,
despite the uncertainty stemming from conjectures in the assumptions
of the theorem.

We believe the time has come for the area to embrace the uncertainty
as well.  Having to rely on unproven assumptions can be rather uncomfortable,
of course, but it is preferable over the alternatives, such as having a
blind conviction, or believing in nothing at all.  We are reminded of
Voltaire's famous dictum made in a theological context:

\smallskip

\begin{center}\begin{minipage}{13.9cm}%
{\emph{``Le doute n'est pas un \'etat bien agr\'eable, mais l'assurance est un \'etat ridicule.''} \\
``Doubt is not a pleasant condition, but certainty is an absurd one,'' Voltaire, 1770.}
\end{minipage}\end{center}

\vskip.75cm

%
%
%




\end{document}